\magnification=\magstephalf
\hsize=150mm
\vsize=210mm

\font\chaptbf=cmbx10 scaled 1400
\font\sc=cmcsc10

\def\chapter#1#2\par{\vfill\eject
  \vglue1in
 \begingroup
 \hyphenpenalty=10000
 \baselineskip=1.5\baselineskip
 \parindent=0pt\leftskip=2pc\rightskip=30mm plus 1fil
    \chaptbf #2\endgroup\vskip10mm}

\def\section#1#2\par{\begingroup \rightskip=0pt plus 1fil \bf #1\endgroup
#2 \bigskip}

\def\citeiteitem#1#2#3\par{
 \leftskip=15mm
 \noindent\llap{\rm #1\ \ }{\sc #2}\ #3\par}

\font\eightrm=cmr8

\font\bbbk = msbm10
\font\eufm = eufm10

\def\Bbbk{\hbox{\bbbk \char"7C}}
\def\Bbb #1{\hbox{\bbbk #1}}
\def\goth #1{\hbox{\eufm #1}}

\def\idd{\mathop {\fam 0 id}\nolimits}

\def\Coeff{\mathop {\fam 0 Coeff}}

\def\o#1{\mathrel{{\circ }_{#1}} }
\def\tag#1{\eqno{(#1)}}


\def\idd{\mathop {\fam 0 id}\nolimits}

\def\Coeff{\mathop {\fam 0 Coeff}}

\def\binom#1#2{\left(\matrix{#1\cr#2\cr}\right)}

\chapter{}{TRIVIAL CONSTRUCTION OF FREE ASSOCIATIVE\hfill
 \lower 10pt\hbox{}\break%
CONFORMAL ALGEBRA}

{\sc Pavel Kolesnikov}
\vskip10mm

{\narrower \eightrm \noindent
This note is to show the effectiveness of the notion of
pseudoalgebra in the theory of conformal algebras.
We adduce very simple construction of free associative
conformal algebra and find its linear basis.
There is no any new result
but we hope that the technique could be useful for further
development of the theory of conformal algebras and pseudoalgebras.
\par}

\vglue 10mm
\section{1\quad INTRODUCTION}

Theory of conformal algebras (see [K1], [K2], [K3])
is a relatively new branch of algebra closely related to mathematical
physics. The general categorial approach in this theory leads to
the notion of pseudotensor category [BD] (also known as multicategory [La]).
Algebras in these categories (so called pseudoalgebras, see [BDK])
allow to get a common presentation of various features of usual and
conformal algebras.

Free associative, commutative and Lie
conformal algebras were investigated in [Ro1], [Ro2] by
using their coefficient algebras.
There were found the bases of
free associative conformal algebra
and its  coefficient algebra (positive component).

Another construction of free associative conformal
algebra was given in [BFK], where it was built without
using coefficient algebra.

In this paper, we adduce another (and quite short)
construction of free associative conformal algebra
and discuss some possibilities of applying some analogous
constructions in the theory of conformal algebras.

This work was supported by RFBR, project 01--01--00630.
Author is very grateful to L.~A.~Bo\-kut, I.~V.~L'vov, V.~N.~Zhelyabin and
E.~I.~Zel'manov for their interest to this work and helpful discussions.

\vglue 10mm
\section{2\quad MAIN DEFINITIONS}

{\sc 2.1. Conformal Algebras}
\smallskip

Let $\Bbbk$ be a field of zero characteristic,
$\Bbbk[D]$ be a polynomial ring on one variable.
A left $\Bbbk[D]$-module $C$ endowed with a family of
$\Bbbk$-bilinear maps
$$
\o{n} : C\otimes C \to C, \quad n\in \{0,1,2,\dots \},
$$
is called {\it conformal algebra} if it satisfies the
following axioms:
$$
 \hbox{\bf locality:}\quad \forall a,b\in C \quad \exists N\ge 0 :
a\o{n} b = 0 \hbox{ for } n\ge N;
                              \tag{1}
$$
$$
 \hbox{\bf sesqui-linearity:}\quad
Da \o{n} b = -n a\o{n-1} b, \quad a\o{n} Db = D(a\o{n} b) +n a\o{n-1} b.
                            \tag{2}
$$
The minimal $N=N(a,b)$ satisfying (1) for fixed $a,b\in C$
defines
{\it locality function} on~$C$.

Conformal algebra $C$ is said to be associative if the
following relations hold:
$$
(a\o{n} b)\o{m} c = \sum\limits_{s\ge 0} (-1)^s
\binom{n}{s} a\o{n-s} (b\o{m+s} c),
\quad  n,m\ge 0,\ a,b,c\in C.
                                            \tag{3}
$$
For detailed explanation of the last notion
one can see [K1], [K3], [Ro1]. In brief,
conformal algebra is associative iff its coefficient algebra
is an associative algebra.

\medskip
{\sc 2.2. Pseudoalgebras}
\smallskip

Let $H$ be a Hopf algebra (see, e.g., [Sw])
with comultiplication $\Delta $, counit $\varepsilon $
and antipode $S$.
We will use the following notation:
instead of $\Delta (h) = \sum\limits_{i} h'_i \otimes h''_i\in H\otimes H$
we will write
$$
\Delta (h) = \sum\limits_{(h)} h_{(1)}\otimes h_{(2)} = h_{(1)}\otimes h_{(2)}
$$
(just eliminating the symbol $\sum$).
Then, because  $\Delta$ is coassociative,
we can denote
$$
(\Delta\otimes \idd)\Delta(h)
=(\idd\otimes \Delta)\Delta(h) = h_{(1)}\otimes h_{(2)} \otimes h_{(3)}
$$
and so on.

The tensor product $H^{\otimes n}$ has a natural structure
of right $H$-module:
$$
(h_1\otimes\dots \otimes h_n)\cdot f =
(h_1\otimes\dots \otimes h_n)\Delta^{[n]} f =
h_1f_{(1)}\otimes \dots \otimes h_nf_{(n)}.
$$

Hopf algebra $H$ is said to be cocommutative if
$h_{(1)}\otimes h_{(2)} = h_{(2)}\otimes h_{(1)}$.
One have to note that for cocommutative Hopf algebras
every permutation of tensor factors in $H^{\otimes n}$
is an endomorphism of right $H$-module defined above.

For our purposes, it would be enough to consider Hopf algebra
$H=\Bbbk[D]$, where
$\Delta (D)= D\otimes 1 + 1\otimes D$,
$\varepsilon (D)=0$, $S(D)=-D$.
Then it is easy to see that
$$
 \Delta(f(D)) = \sum\limits_{s\ge 0} D^{(s)} \otimes {d^s f\over dD^s},
$$
where $D^{(s)}={1\over s!}D^s$. Hereinafter we will use
$x^{(s)}$ for ${1\over s!}x^s$.

An associative algebra $A$ endowed with
coaction $\Delta_A$ of a Hopf algebra $H$ is called {\it $H$-comodule algebra}
if the coaction is a homomorphism of algebras.
More carefully, $A$ is an $H$-comodule algebra if
the homomorphism of algebras
$\Delta_A : A \to H\otimes A$ satisfies
$(\Delta\otimes \idd_A)\Delta_A = (\idd_H\otimes \Delta_A)\Delta_A$
and
$(\varepsilon\otimes \idd_A)\Delta_A(a) = 1\otimes a$.

\medskip
{\sc Definition {\rm [BD], [BDK]}.}
A left $H$-module $P$ endowed with $H$-bilinear map
$$
 * : P\otimes P \to H\otimes H \otimes _H P
                                                    \tag{4}
$$
is called a {\it pseudoalgebra} (or $H$-pseudoalgebra).
\medskip

It is suitable to expand $*$ by the following way:
consider
$$
* : (H^{\otimes n}\otimes _H P)\otimes (H^{\otimes m}\otimes _H P)
\to H^{\otimes n+m} \otimes _H P
$$
via
$$
 (F\otimes_H a)*(G\otimes _H b) =
(F\otimes G\otimes_H 1)(\Delta^{[n]}\otimes \Delta^{[m]}\otimes_H \idd)(a*b).
                        \tag{5}
$$

{\bf Proposition~1 {\rm [BDK]}.}
{\sl
Every conformal algebra $C$
is an $H=\Bbbk[D]$-pseudoalgebra with the
pseudoproduct
$$
 a*b = \sum\limits_{s\ge 0} (-D)^{(s)}\otimes 1 \otimes_H (a\o{n} b).
                                                     \tag{6}
$$

Every $H=\Bbbk[D]$-pseudoalgebra is a conformal algebra
with the $\o{n}$-products defined by~{\rm (6)}.
}
\medskip

{\sc Proof.}
The first statement could be easily deduced from the axioms
(1),~(2) of conformal algebra.
To check the second statement, one should note that
every element
$f\otimes g\in \Bbbk[D]\otimes \Bbbk[D]$
could be uniquely presented as
$$
f\otimes g = \sum\limits_{t\ge 0} ((-D)^{(t)}\otimes 1)\Delta (h_t),
\quad h_t\in \Bbbk[D].
$$
Namely (see, e.g., [BDK] for more general fact),
$$
(f\otimes g) = \sum\limits_{s\ge 0} ((-D)^sf\otimes 1)
\Delta\left ({d^s g \over dD^s}\right).
$$
This relation has a simple interpretation in terms of
Hopf algebras:
$f\otimes g = (fS(g_{(1)})\otimes 1)\Delta(g_{(2)})$.

Then, locality axiom (1) follows directly from the definition
of pseudoproduct~(4); sesqui-linearity (2) follows from
$H$-bilinearity of pseudoproduct. \quad $\diamond $
\medskip

An $H$-pseudoalgebra $P$ is said to be associative if
$$
a*(b*c) = (a*b)*c ,
                                                                     \tag{7}
$$
see~(5). It is clear, an associative
conformal algebra is just the same as associative
$\Bbbk[D]$-pseudoalgeb\-ra [BDK].

\medskip
{\bf Proposition 2 {\rm [Ko1]}.}
{\sl
Let $H$ be a cocommutative Hopf algebra and $A$ be an
$H$-comodule algebra.
Then free $H$-module $\goth A = H\otimes A$
with a pseudoproduct  defined by
$$
 (h\otimes a)*(g\otimes b) = (hb_{(1)}\otimes g)\otimes _H (1\otimes ab_{(2)})
                                                                     \tag{8}
$$
or
$$
 (h\otimes a)*(g\otimes b) = (h\otimes gS(a_{(1)}))\otimes _H (1\otimes a_{(2)}b)
                                                                     \tag{9}
$$
is an associative pseudoalgebra.
}

\medskip
{\sc Proof.}
Because of $H$-bilinearity of (8) and (9), it is sufficient to show (7)
for elements of $1\otimes A$.

Let us prove the proposition for pseudoproduct~(8), for example.
Straightforward calculation shows (see (5)) that
$$
((1\otimes a)*(1\otimes b))*(1\otimes c)
= (1\otimes S(a_{(1)})\otimes S(b_{(1)})S(a_{(2)}))\otimes_H (1\otimes a_{(3)}b_{(2)}c)
= (1\otimes a)*((1\otimes b)*(1\otimes c)).  \quad \diamond
$$

By the same arguments, we have

\medskip
{\bf Proposition 3 {\rm [Ko1]}.}
{\sl
Let $H$ be a commutative Hopf algebra and $A$ be an
$H$-comodule algebra.
Then $\goth A = H\otimes A$
with a pseudoproduct  defined by
$$
 (h\otimes a)*(g\otimes b) = (h\otimes ga_{(1)})\otimes _H (1\otimes a_{(2)}b )
                                                                     \tag{10}
$$
is an associative pseudoalgebra. \quad $\diamond $
}

\medskip
{\bf Proposition 4.}
{\sl
Let $H$ be a commutative and cocommutative Hopf algebra and $A$ be an
$H$-comodule algebra.
Then $\goth A = H\otimes A$
with a pseudoproduct  defined by
$$
 (h\otimes a)*(g\otimes b) = (hS(b_{(1)})\otimes g )\otimes _H (1\otimes ab_{(2)})
                                                                     \tag{11}
$$
is an associative pseudoalgebra. \quad $\diamond $
}
\medskip

{\sc Example 1.}
Consider an associative algebra $A$ and
$\Delta_A : a \mapsto 1\otimes a$. This coaction
turns $H\otimes A$ to be a {\it current pseudoalgebra} (see, e.~g., [K1]).
\medskip

{\sc Example 2.}
Let $H=\Bbbk[D]$ and $A=\Bbbk[v]$,
$\Delta_A (v) = D\otimes 1 + 1\otimes v$.
Then $H\otimes A$ endowed with any of the pseudoproducts (8)--(11)
is known as {\it Weyl conformal algebra}. Its adjoint conformal Lie algebra,
i.e., the same module endowed with ``pseudocommutator''
$$
[a*b] = a*b - (\sigma_{12}\otimes _H \idd)(b*a),
\quad
\sigma_{12}(h_1\otimes h_2) = h_2\otimes h_1,
$$
contains so called {\it Virasoro conformal algebra}.
\medskip

\vglue 10mm
\section{3\quad FREE ASSOCIATIVE CONFORMAL ALGEBRA}

The purpose of this note is to show ``universality'' of the construction described
in propositions~2,~3.
In this section, $\Bbbk$ be a field of zero characteristic
and $H=\Bbbk[D]$.

It is clear there are no universal objects in the class of all
associative conformal algebras (with a fixed set~$B$ of generators).
But if we fix some locality function
$N : B\times B \to \Bbb Z_+$
then it becomes possible
to build free algebra $CF(B,N)$ in the class of associative
conformal algebra generated by $B$ with locality function~$N$ (see [Ro1]).
Direct construction of such a conformal algebra is given in [BFK] (for
$N=\hbox{const}$).
In these papers,
it was proved that conformal monomials
$$
a_1\o{n_1}(a_2\o{n_2} \dots \o{n_{k-1}} (a_k\o{n_k} D^s a_{k+1})\dots ),
\quad
a_i\in B,
\
0 \le n_i < N(a_i, a_{i+1}),
\
s\ge 0,
\ k\ge 0,
                                                  \tag{12}
$$
form a linear basis of
$CF(B, N)$.
In this note, we adduce another direct construction of this algebra:
for $N(a,b)=N(b)$, i.e., locality function depends only on its second
argument.

Consider an arbitrary set of symbols (generators) $B$;
let $F(B)= \Bbbk \langle B\cup\{v\} \rangle$ be the
free associative algebra generated by $B$ and an additional element~$v$
(``Virasoro element'').

Define on $F(B)$ the following comodule algebra structure:
$$
\displaylines{
\qquad\qquad\qquad\qquad
\llap{$\Delta_F :\ $} F(B) \to H\otimes F(B), \hfill\cr
\qquad\qquad\qquad\qquad
a\mapsto 1\otimes a,  \quad a\in B,\hfill\cr
\qquad\qquad\qquad\qquad
v\mapsto D\otimes 1 + 1\otimes v.  \hfill\cr
}
$$
Since $F(B)$ is free, it is sufficient to define $\Delta_F$ only on
the generators of this algebra.

Denote by
$\goth F(B)$ the pseudoalgebra $H\otimes F(B)$
endowed with the pseudoproduct~(8).

\medskip
{\bf Theorem~1.}
{\sl
Let us fix positive function
$n(a)\in \Bbb Z_+$, $a\in B$.
Then the pseudoalgebra
$\goth F_n(B)$
generated by
$\{1\otimes v^{(n(a)-1)}a \mid a\in B\}$
is isomorphic to the free associative
conformal algebra generated by $B$ with
locality function $N(a,b)=n(b)$.
}
\medskip

{\sc Proof.}
Let $C$ be an associative conformal algebra
generated by the set $B$ with locality function $N$.
Conformal monomials
$$
D^s(a_1\o{n_1}(a_2\o{n_2} \dots \o{n_{k-1}} (a_k\o{n_k} a_{k+1})\dots )),
\quad k\ge 0,
                                                  \tag{13}
$$
are called {\it normal words} if
$$
a_i\in B,
\quad
0 \le n_i < N(a_i, a_{i+1}),
\quad
s\ge 0.
                                                  \tag{14}
$$
Normal word is {\it $D$-free} if $s=0$.

\medskip
{\bf Lemma 1.}
{\sl
Let $w$ be a $D$-free normal word and $a_0\in B$.
Then for every $n\ge 0$ the element
$a_0\o{n} w\in C$
could be represented as a $\Bbbk$-linear combination of
$D$-free normal words starting with symbol~$a_0$.
}

\medskip
{\sc Proof.}
Let
$$
w=a_1\o{n_1}(a_2\o{n_2} \dots \o{n_{k-1}} (a_k\o{n_k} a_{k+1})\dots ),
k\ge 0.
$$
For $k=0$ the statement is obvious.
If $k>0$ and $n<N(a_0, a_1)$ then $a_0\o{n} w$ is a~normal word.
So it remains to consider
$k>0$, $n\ge N(a_0, a_1)$.
Denote $w=a_1\o{n_1} w_1$.

It is easy to note that conformal associativity condition~(3) is equivalent
to the following one:
$$
a\o{n} (b\o{m} c) = \sum\limits_{s\ge 0} \binom{n}{s} (a\o{n-s} b)\o{m+s} c.
$$
Hence,
$$
\displaylines{
\qquad
a_0\o{n} (a_1\o{n_1} w_1)
 = \sum\limits_{s\ge n-N(a_0,a_1)+1}
   \binom{n}{s} (a_0\o{n-s} a_1)\o{n_1+s} v_1
   \hfill \cr
 \hfill
 = \sum\limits_{s\le N(a_0,a_1)-1}
    \alpha_s (a_0\o{s} (a_1 \o{n+n_1-s} v_1)),
    \quad
    \alpha_s \in \Bbbk. \hfil\qquad \hbox{\rm(15)}\cr
}
$$
Since
$w_1$ is shorter than~$w$ then we can assume
$a_1 \o{n + n_1 - s} w_1$
to be representable as a linear combination of normal words
starting with~$a_1$. Namely,
 let
$$
a_1 \o{n + n_1 - s} w_1 = \sum\limits_{t} \beta_{s, t} w_{s,t},
$$
where $w_{s,t}$ are $D$-free normal words starting with $a_1$. Then
$$
a_0\o{n} w = \sum\limits_{s\le N(a_0,a_1)-1} \sum\limits_t
  \beta_{s,t}  \alpha_s (a_0\o{s} w_{s,t}),
$$
where
$a_0\o{s} w_{s,t}$ are normal words. \quad $\diamond $

\medskip
{\bf Lemma 2 {\rm [BFK]}.}
{\sl
Every element $f\in C$ could be represented as
a~linear combination of normal words.
}
\medskip

{\sc Proof.}
Let us define {\it conformal word} by the following way:

i) $a\in B$ is a conformal word;

ii) if $u$ and $v$ are conformal words then $u\o{n} v$ is also conformal word for every
$n\ge 0$.

It is clear, every element $f\in C$ could be represented as a
$\Bbbk[D]$-linear combination of conformal words (see~(2)).
Using associativity~(3) one can rewrite every conformal word
as a linear combination of right-normed words~(13).
Now it is sufficient to apply lemma~1
in order to represent every right-normed word in normal form~(13), (14).
\quad $\diamond $
\medskip

Consider a linear order $\le $ on $B$ and extend it to
$B\cup \{v\}$ by assuming $a< v$ for every $a\in B$.
Then basic monomials of
$F(B)$
could be linearly ordered by {\it deg-lex rule}:
$$
x_1 \dots x_k \le y_1\dots y_m \iff
 k<m \ \hbox{or}\ k=m \ \hbox{and}\
  x_1 \dots x_k \ \hbox{is lexicographically less than}\  y_1\dots y_m.
$$
For  every $0\ne f\in F(B)$ we can determine its principal monomial
$\bar f$:
$f= \alpha \bar f + \sum_s \alpha_s u_s$, $\alpha\ne 0$, $u_s<\bar f$.
 Also, one can define lower monomial $\hat f$ for every
$0\ne f\in F(B)$.

The following lemma is obvious.

\medskip
{\bf Lemma 3.}
{\sl
Let $u_1,u_2,u_3$ be some monomials in $F(B)$.

 {\rm 1.} If $u_1\le u_2$ then  $u_3u_1\le u_3u_2$
and $u_1u_3 \le u_2u_3$.

 {\rm 2.} If $u_1\le u_2$, $f_i={du_i\over dv}$, $i=1,2$,
 then $\hat f_1\le \hat f_2$.

For both of the statements strict inequalities also hold.
\quad $\diamond $
}
\medskip

It follows from lemma~3 that for every polynomials
$f,g\in F(B)$, $\hat f \le \hat g$,
we have
$\hat f' \le hat g'$, where $f'={df\over dv}$, $g'={dg\over dv}$.

For fixed function $n:B\to \Bbb N$,
we identify $1\otimes v^{(n(a)-1)}a\in \goth F_n(B)$
with $a\in B$.
Let us calculate monomial (13) under the conditions~(14) in
$\goth F_n(B)$.

First, we note that for every $f\in F(B)$
$$
(1\otimes v^{(n(a)-1)}a)*(1\otimes f) =
\sum\limits_{s\ge 0} (D^{(s)}\otimes 1)\otimes_H
 \left (1\otimes v^{(n(a)-1)}a{d^s f \over dv^s}\right).
$$
Hence,
$$
(1\otimes v^{(n(a)-1)}a)\o{m} (1\otimes f) = (-1)^m
\left (1\otimes v^{(n(a)-1)}a \left ({d^m f \over dv^m}\right)\right).
                                                      \tag{16}
$$
Let
$u=a_1\o{n_1}(a_2\o{n_2} \dots \o{n_{k-1}} (a_k\o{n_k} a_{k+1})\dots )$
be a $D$-free normal monomial in $\goth F_n(B)$.
It is easy to see from (16) that
$$
 u=1\otimes f(a_1,\dots a_{k+1}; n_1,\dots , n_k).
$$
 Lemma 3 implies that
$$
\displaylines{
\quad \quad \quad
 \hat f(a_1,\dots a_{k+1}; n_1,\dots , n_k) \hfill \cr
\hfill
=(-1)^{n_1+\dots n_k}
v^{(n(a_1)-1)}a_1 v^{(n(a_2)-n_1-1)}a_2 \dots
a_k  v^{(n(a_{k+1})-n_k-1)}a_{k+1}\hfil
\quad \quad \quad
                                                    \llap{(17)}\cr
}
$$
Here $n(a_{i+1})-n_i-1\ge 0$ since the conditions~(14) hold.

Monomials of the form~(17) are linearly independent
in free associative algebra $F(B)$. Hence,
polynomials $f(a_1,\dots a_{k+1}; n_1,\dots , n_k)$
are also linearly independent. Therefore, we have constructed
associative conformal algebra
$\goth F_n(B)$
generated by~$B$
with locality $N(a,b)=n(b)$ such that
normal words are linearly independent in this algebra.
Lemma 2 implies that every associative conformal algebra
generated by $B$ with locality function less or equal $N(a,b)=n(b)$
is a homomorphic image of $\goth F_n(B)$.
This accomplish the proof of theorem~1. \quad $\diamond $
\eject

\vglue 10mm
\section{4\quad ANOTHER COMODULAR CONSTRUCTION}

In this section, we consider another constructions of
pseudoalgebra over commutative and cocommutative Hopf algebras.
Over algebraically closed field $\Bbbk$ of zero characteristic, the only examples of
such Hopf algebras are
$H= \Bbbk[X]\otimes \Bbbk[\Gamma]$, where $X$ is a set of
commuting generators and $\Gamma $ is an abelian group
(see, e.g., [Sw]). In particular this  construction will
be applicable to conformal algebras ($H=\Bbbk[D]$).

Let $A$ be an algebra (not necessarily associative), $H$
be a Hopf algebra.
Then the homomorphism of algebras
$$
 \Delta_A : A \to H\otimes A, \quad  a\mapsto a_{(1)}\otimes a_{(2)},
$$
we call by {\it coaction} of $H$ on $A$
if
$$
(\Delta \otimes \idd_A)  \Delta_A(a) =
(\idd_H\otimes \Delta_A) \Delta_A(a),
\quad
\varepsilon (a_{(1)}) a_{(2)} = a.
$$

\medskip
{\bf Proposition~5 {\rm [Ko2]}.}
{\sl
 Let $C$ be a conformal algebra. Then the following
conditions are equivalent.

{\rm 1.} Coefficient algebra  $\Coeff C$
satisfies poly-linear homogeneous identities of the form
$$
 \sum\limits_{\sigma \in S_n} t_\sigma (a_{1\sigma },\dots a_{n\sigma })=0,
                                                    \tag{18}
$$
where
$t_\sigma(y_1,\dots ,y_n)$
is a linear combination of non-associative words obtained
from  $y_1\dots y_n$
by some bracketing.

{\rm 2.} Considered as a pseudoalgebra,
$C$
satisfies
``pseudo''-identities of the form
$$
\sum\limits_{\sigma \in S_n}
(\sigma \otimes_H \idd)t^*_\sigma (a_{1\sigma },\dots a_{n\sigma })=0,
                                                    \tag{19}
$$
where
 $t^*$
 means the same term
 $t$ with operation~$*$ instead of usual multiplication. \quad $\diamond$
}

\medskip
{\bf Theorem~2.}
{\sl
Let $H$ be a commutative and cocommutative
Hopf algebra, $A$ be an algebra endowed with coaction of $H$.
Then $H\otimes A$ with the pseudoproduct defined by
$$
 (h\otimes a)*(g\otimes b) =
 (hb_{(1)}\otimes ga_{(1)})\otimes _H(1\otimes a_{(2)}b_{(2)})
                                                    \tag{20}
$$
is an
$H$-pseudoalgebra.

If $A$ satisfies identity~(18) then the pseudoalgebra
$H\otimes A$ with the pseudoproduct~{\rm (20)}
satisfies~{\rm (19)}.
}
\medskip

{\sc Proof.}
Let $t(a_1,\dots ,a_n)$ be a non-associative word
obtained from $a_1\dots a_n$ by some bracketing.
It is sufficient to prove that
$$
 t^*(a_1,\dots ,a_n)=
\left (\bigoplus\limits_{k=1}^n
a_{1(k-1)}\dots a_{k-1(k-1)} a_{k+1(k)} \dots  a_{n(k)}
\right)
\otimes_H (1\otimes t(a_{1(n)},\dots ,a_{n(n)})).
                                                    \tag{21}
$$
It could be easily done by using induction on~$n$. \quad $\diamond $
\medskip

{\bf Corollary.}
{\sl
Let $F[B]$ be a free associative commutative
algebra on generators $B\cup\{v\}$.
Define coaction of $H=\Bbbk[D]$ on $F[B]$:
$$
\displaylines{
\qquad\qquad\qquad\qquad
\llap{$\Delta_F :\ $} F[B] \to H\otimes F[B], \hfill\cr
\qquad\qquad\qquad\qquad
a\mapsto 1\otimes a,  \quad a\in B,\hfill\cr
\qquad\qquad\qquad\qquad
v\mapsto D\otimes 1 + 1\otimes v.  \hfill\cr
}
$$
Then the pseudoalgebra
$\goth F[B]$ with the pseudoproduct~{\rm (20)} is
an associative and commutative conformal algebra.
}
\medskip

{\bf Question 1.}
How to construct free associative commutative conformal algebra
using the construction described in Theorem~2?
\medskip

{\bf Question 2.}
How to construct free Lie conformal algebra via the same way?

\vglue 10mm
\section{REFERENCES}

\citeiteitem{[BDK]}{Bakalov~B., D'Andrea~A., Kac~V.~G.}Theory of finite pseudoalgebras, Adv. Math., {\bf 162}, N.~1 (2001).

\citeiteitem{[BD]}{Beilinson~A.~A., Drinfeld~V.~G.}Chiral algebras, Preprint.

\citeiteitem{[BFK]}{Bokut~L.~A., Fong~Y., Ke~W.-F.}Free associative conformal algebras, Proc. of the 2nd Tainan-Moscow Algebra and Combinatorics Workshop, Tainan 1997, 13--25. Springer-Verlag, Hong Kong, 2000.

\citeiteitem{[K1]}{Kac~V.~G.}Vertex algebras for beginners, Univ. Lect. Series {\bf 10}, AMS, 1996.

\citeiteitem{[K2]}{Kac~V.~G.}The idea of locality. In: H.-D. Doebner et al. (eds.) Physical applications and mathematical aspects of geometry, groups and algebras, Singapore: World Scientific (1997), 16--32.

\citeiteitem{[K3]}{Kac~V.~G.}Formal distribution algebras and conformal algebras, XII-th International Cong\-ress in Mathematical Physics (ICMP'97) (Brisbane), Internat. Press: Cambridge, MA (1999) 80--97.

\citeiteitem{[Ko1]}{Kolesnikov P.~S.}On universal representations of Lie conformal superalgebras. Novosibirsk, 2002 (Preprint~/ SB RAS, Sobolev Inst. Math., N~102).

\citeiteitem{[Ko2]}{Kolesnikov P.~S.}Simple Jordan pseudoalgebras. Novosibirsk, 2002 (Preprint~/ SB RAS, Sobolev Inst. Math., N~103). ArXive math.QA/0210264.

\citeiteitem{[La]}{Lambek~J.}Deductive systems and categories. II. Standard constructions and closed categories, Lecture Notes Math., {\bf 86}, 76--122, Springer-Verl., Berlin, 1969.

\citeiteitem{[Ro1]}{Roitman~M.}On free conformal and vertex algebras, J. Algebra, {\bf 217} (1999), 496--527.

\citeiteitem{[Ro2]}{Roitman~M.}Universal enveloping conformal algebras, Sel. Math., New Ser., {\bf 6}, N.~3 (2000), 319--345.

\citeiteitem{[Sw]}{Sweedler~M.}Hopf Algebras. Math. Lecture Note Series, New York: Benjamim, 1969.
\medskip \medskip

\leftline{Sobolev Institute of Mathematics}
\leftline{E-mail: pavelsk@math.nsc.ru}

\end